\documentclass{numapde-preprint}

\usepackage{numapde-semantic}
\usepackage{numapde-style}
\usepackage{numapde-local}

\hypersetup{
	pdftitle={An Extension of the Strain Transfer Principle for Fiber Reinforced Materials},
	pdfauthor={Felix Ospald, Kai Bergermann, Roland Herzog},
	pdfkeywords={strain transfer principle, strain measurement, fiber reinforced materials, fiber Bragg sensors, fiber optical strain sensors}
}

\title{An Extension of the Strain Transfer Principle for Fiber Reinforced Materials\thanks{This work is part of a measure which is co-financed by tax revenue based on the budget approved by the members of the Saxon state parliament. Financial support is gratefully acknowledged.}}
\subtitle{}
\shorttitle{Extended Strain Transfer Principle}

\author{Felix Ospald\thanks{Technische Universität Chemnitz, Faculty of Mathematics, 09107 Chemnitz, Germany (\email{felix.ospald@mathematik.tu-chemnitz.de}, \url{https://www.tu-chemnitz.de/mathematik/part_dgl/people/ospald}, \orcid{0000-0001-8372-9179}, \email{kai.bergermann@mathematik.tu-chemnitz.de}, \url{https://www.tu-chemnitz.de/mathematik/part_dgl/people/bergermann}, \orcid{0000-0002-2259-1839}, \email{roland.herzog@mathematik.tu-chemnitz.de}, \url{https://www.tu-chemnitz.de/mathematik/part_dgl/people/herzog}, \orcid{0000-0003-2164-6575}).}
\and
Kai Bergermann\footnotemark[2]
\and
Roland Herzog\footnotemark[2]}
\shortauthor{F. Ospald, K. Bergermann and R. Herzog}

\dedication{}

\IfFileExists{numapde-uncolorizeHyperrefIfChanges.sty}{\RequirePackage{numapde-uncolorizeHyperrefIfChanges}}{}

\begin{document}
\maketitle

\begin{abstract}
Fiber optical strain sensors are used to measure the strain at a particular sensor position inside the fiber. 
In order to deduce the strain in the surrounding matrix material, one can employ the strain transfer principle.
Its application is based on the assumption that the presence of the fiber does not impede the deformation of the matrix material in fiber direction.
In fact, the strain transfer principle implies that the strain in fiber direction inside the fiber carries over verbatim to the strain inside the matrix material. 
For a comparatively soft matrix material, however, this underlying assumption may not be valid.
To overcome this drawback, we propose to superimpose the matrix material with a one-dimensional model of the fiber, which takes into account its elastic properties. 
The finite element solution of this model yields a more accurate prediction of the strain inside the fiber in fiber direction at low computational costs.\end{abstract}




\section{Introduction}
\label{sec:introduction}

Fiber optical sensors, such as fiber Bragg gratings embedded into a surrounding matrix material, are often used to measure the strain at the sensor position inside the fiber. 
Such measurements can be used, for instance, to infer the magnitude of residual stresses in the matrix material.
This can be achieved through an appropriate inverse problem, based on a forward deformation simulation.
However, the numerical simulation of the stresses and strains inside parts with an embedded measurement fiber under mechanical loads is challenging due to the difference in typical length scales between the fiber diameter and part geometry.
A potential way out is to simulate the matrix material in the absence of the fiber, and to incorporate the effect of the latter only a~posteriori.

The first attempt to analytically model the stress transfer from a uniaxially loaded matrix material to an embedded fiber was made by \cite{Cox:1952:1}, leading to the emergence of the research field referred to as shear-lag theory today (see for instance \cite{Nayfeh:1977:1,McCartney:1992:1,Nairn:1997:1}), with various applications to fiber optical sensors described, \eg, in \cite{LiZhouLiangLi:2007:1,LiZhouRenLi:2009:1,ZhouLiuHuangWangJianpingHuangJinping:2012:1}. 
Here, the uniaxial stress in fiber direction is related to the shear stress at the fiber matrix interface, which is recognized as the dominating mechanism of stress transfer from the matrix to the fiber material. 

In order to deduce the full strain state inside the material surrounding the fiber instead of only two stress components, one can alternatively employ the strain transfer principle (STP) described in \cite{Lekhnitskii:1981:1,KollarVanSteenkiste:1998:1}.
The STP postulates a linear relationship between the strain tensor inside the sensor and the strain tensor of the far field of the surrounding matrix material (\ie, as though there was no fiber present).
This linear relationship can be expressed analytically in the form of the strain transfer tensor and it is valid for orthotropic matrix as well as orthotropic fiber materials; see \cite{KollarVanSteenkiste:1998:1}. 
An extension of this model to coated fibers and temperature differences between the matrix and the fiber is also available in \cite{VanSteenkisteKollar:1998:1}. 
The predictions of the analytical STP have been confirmed by various experimental works, for instance \cite{LuyckxVoetWaeleDegrieck:2010:1,VoetLuyckxWaeleDegrieck:2010:1,LammensLuyckxVoetVanPaepegemDegrieck:2015:1}. 

The STP yields particularly good results when the material properties of the matrix and the fiber are similar, or when the fiber material is softer than the matrix. 
In these cases the fiber does not restrain the deformations of the matrix material in fiber direction under a certain load. 
The strain in fiber direction in the matrix material is transferred verbatim to the strain inside the fiber in fiber direction and vice versa.
In case the matrix material is softer than the fiber, however, the fiber itself may restrain deformations of the entire part/matrix material, and the strain in fiber direction no longer carries over verbatim from the matrix material.
The magnitude of this effect also depends on the fiber diameter, part dimensions as well as load conditions under consideration.

To overcome this drawback one has to take into consideration the entire geometry of the part, the embedded fiber as well as the load conditions. 
However, since the fiber diameter is usually small compared to the part dimensions, a fully resolved finite element (FE) model is often impractical.
Instead, we propose an extension of the STP.
Our method combines the practical benefits of the STP with the improved accuracy of a fully resolved finite element model.
In particular, we can continue to simulate the deformation of the matrix material in the absence of the fiber. 

We refer to our proposed approach as the extended STP.
We apply the classical STP to deduce all strain components except the strain in fiber direction from a FE model of the matrix material without the fiber.
By contrast, the strain in fiber direction is derived from the solution of a modified FE model.
The latter is obtained by superimposing the elastic properties of the bulk matrix material with a one-dimensional model of the fiber.
This does not require the fiber to be resolved in the computational mesh.

The paper is structured as follows. \Cref{sec:linear_elasticity} states the linear elasticity problem for a matrix material part with an embedded fiber. \Cref{sec:strain_transfer_principle} introduces the strain transfer principle and recalls the results from the existing analytical theory.
\Cref{sec:strain_transfer_principle_extensions} presents our extension to this theory. 
In \cref{sec:variational_form} we derive the variational form of our elasticity model, and \Cref{sec:experiments} presents detailed numerical results of our extended STP in comparison with the original STP and with fully resolved finite element computations.

\textbf{Nomenclature:} 
We denote by $A \dprod B$ the double contraction of a rank-$4$ tensor $A$ with a matrix $B$, \ie, $(A \dprod B)_{ij} = \sum_k \sum_\ell A_{ijk\ell} B_{k\ell}$. 
We also use $A \dprod B$ for the double contraction of two matrices, \ie, $A \dprod B = \sum_i \sum_j A_{ij} B_{ij}$.
Moreover, $\I$ denotes the identity on rank-$4$ tensors, \ie $\I \dprod A = A$ holds for any matrix~$A$ of appropriate dimensions.
$A \otimes B$ denotes the outer product between matrices $A$ and $B$, \ie $(A \otimes B)_{ijk\ell} = A_{ij} B_{k\ell}$.
Finally $a \cdot b$ denotes the usual dot product between vectors $a$ and $b$, \ie, $a \cdot b = \sum_i a_i b_i$.

\section{Linear Elasticity with Embedded Fiber}\label{sec:linear_elasticity}

Let $\Omega \subset \R^d$ denote the domain ($d=3$) occupied by the part under consideration.
Furthermore, let $\C \colon \Omega \to \Lmap(\Sym(d))$ denote the stiffness tensor field on $\Omega$, \ie for every material point~$x \in \Omega$, $\C(x)$ is a linear map between symmetric $d \times d$ strain matrices and symmetric $d \times d$ stress matrices. The function $f \colon \Omega \to \R^d$ denotes a force density field, \eg, due to gravity. Let $\Uad$ denote the set of all admissible displacements of the part which satisfy the given boundary conditions.

The elastic deformation energy of a displacement field $u \in \Uad$ on the domain without embedded fiber is given by 
\begin{equation}
	\label{eq:matrix_energy}
	\Ematrix(u) \coloneqq \frac{1}{2} \int_{\Omega} \varepsilon(u) \dprod \C \dprod \varepsilon(u) \dx - \int_{\Omega}  f \cdot u \dx,
\end{equation}
where $\varepsilon(u)$ denotes the symmetric displacement gradient
\begin{equation}
	\varepsilon(u) = \frac{1}{2} \left( \nabla u + (\nabla u)^\transp \right).
\end{equation}
The equilibrium solution, which solves the minimization problem 
\begin{equation}
	\label{eq:min-problem_without_fiber}
	\min_{u \in \Uad} \Ematrix(u),
\end{equation}
will be denoted by $\unf$, where the subscript denotes the absence of the fiber.
We refer the reader to \cite{Braess:2007:1} for an account of the mathematical theory.

Let $\gamma \colon [0, L] \to \Omega$ denote the arc-length parameterization of a curve which models the center of the fiber.
Consequently, $L$ denotes the total length of the fiber.
Let us assume that each point of the fiber $\gamma(t)$ is tied to the corresponding point in the domain $\Omega$, \ie, we do not consider slip between fiber and matrix material.
Given a deformation field $u \in \Uad$ on the domain, the energy of a one-dimensional fiber generally consists of three parts,
\begin{equation}
	\Efiber(u) \coloneqq \Estretch(u) + \Ebend(u) + \Etwist(u),
\end{equation}
modelling the stretching energy $\Estretch(u)$, the bending energy $\Ebend(u)$ and the twisting energy $\Etwist(u)$, respectively.
Following \cite{Spencer:1984:1,SpencerSoldatos:2007:1}, the stretching energy is given by
\begin{equation}
	\Estretch(u) \coloneqq \frac{1}{2} \int_{\gamma} E A \varepsilon_\gamma(u)^2 \ds,
\end{equation}
where $E$ denotes the effective elastic modulus of the fiber material defined below, $A$ is the cross-sectional area and $\varepsilon_\gamma(u)$ is the strain in fiber direction.
Similarly, the bending energy is given by
\begin{equation}
		\Ebend(u) \coloneqq \frac{1}{2} \int_{\gamma} E I \kappa_\gamma(u)^2 \ds,
\end{equation}
where $I$ denotes the cross-sectional moment of inertia and $\kappa_\gamma$ denotes the curvature of $\gamma$ (\ie, the derivative of the bending angle).
For the twisting energy we have 
\begin{equation}
	\Etwist(u) \coloneqq \frac{1}{2} \int_{\gamma} G I_T \delta_\gamma(u)^2 \ds,
\end{equation}
where $G$ denotes the shear modulus, $I_T$ is the torsion constant for the section and $ \delta_\gamma(u)$ denotes the derivative of the torsion angle of $\gamma$.
We consider only homogeneous fibers for which $E$, $A$, $I$, $G$ and $I_T$ are constant.
However, the curvature~$\kappa_\gamma$ may vary along the fiber.

Usually, the radius of curvature of the fiber, $1/\kappa_\gamma(u)$, is much larger than the radius~$R$ of the fiber itself. 
Similarly, the length over which the fiber twists by a full turn $2\pi/\delta_\gamma(u)$ is usually much larger than $R$. 
Furthermore, we have $A \propto R^2$, $I \propto R^4$ and $I_T \propto R^4$, from where we conclude $I \kappa_\gamma(u)^2 \ll R^2$ and $I_T \delta_\gamma(u)^2 \ll R^2$.
As a consequence, the bending and twisting energy terms can be neglected compared to the stretching term for most applications involving only small deformations of the domain $\Omega$.

Therefore, we neglect the bending and twisting energies and consider 
\begin{equation}\label{eq:E_fiber}
	\Efiber(u) = \Estretch(u) = \frac{1}{2} \int_{\gamma} E A \varepsilon_\gamma(u)^2 \ds.
\end{equation}
The equilibrium solution, which solves the minimization problem 
\begin{equation}
	\label{eq:min-problem}
	\min_{u \in \Uad} \Ematrix(u) + \Efiber(u),
\end{equation}
will be denoted by $\uf$, where the subscript denotes the presence of the fiber in the deformation energy.

We come back to the definition of the effective stretching Young's modulus of the fiber material.
In order to compensate for the existing matrix material in the volume occupied by the fiber, $E$ is defined as
\begin{equation}
	E = \max \{0, E_f - E_m \},
\end{equation}
where $E_f$ is the stretching Young's modulus of the fiber and $E_m$ is the stretching Young's modulus of the matrix in the local fiber direction $v = \gammadot(s)$.
$E_m$ is calculated from the matrix compliance as 
\begin{equation*}
	E_m = (\projgamma \dprod \C^{-1} \dprod \projgamma)^{-1},
\end{equation*}
where $\projgamma = v v^\transp$. For the case $E_f \le E_m$ which indicates soft fiber material, this results in $\uf = \unf$. 
Otherwise, $\uf$ and $\unf$ are different.

The fiber cross-section is assumed to be of circular shape, such that
\begin{equation*}
	A = \pi R^2
\end{equation*}
holds, and the strain in fiber direction $\varepsilon_\gamma(u)$ is given by
\begin{equation*}
	\varepsilon_\gamma(u) = \projgamma \dprod \varepsilon(u).
\end{equation*}
These relations allow us to evaluate and minimize the total deformation energy in \eqref{eq:min-problem}.
Notice that this does not require to resolve the fiber in the computational mesh.

\section{Strain Transfer Principle}
\label{sec:strain_transfer_principle}

The STP states that there exists a linear relationship between the strain at the center of the fiber and the strain inside the matrix material.
In other words, there exists $\Tnf \in \Lmap(\Sym(d))$ such that 
\begin{equation}
	\label{eq:strain_transfer_principle}
	\varepsilon(\uf) \approx \Tnf \dprod \varepsilon(\unf) \quad \text{along } \gamma.
\end{equation}
$\Tnf$ depends only on the geometry of the fiber described by its radius~$R$ and the material parameters of the matrix and fiber and possibly the fiber orientation if any of the materials are anisotropic.
\Cref{eq:strain_transfer_principle} is exact only if $\unf$ and $\uf$ are identical outside of the fiber. 
In general, this is approximately fulfilled if the fiber has only little influence on the overall deformation $\uf$. 
This is for instance the case if $E_f \ll E_m$, \ie, the fiber is a relatively soft inclusion in the matrix material.

\paragraph{Analytical representation}
\label{sec:analytical_stp}

There exists an analytical representation of the strain transfer principle by \cite{KollarVanSteenkiste:1998:1}, which considers a homogeneous, uncoated orthotropic fiber with elliptic cross-section embedded into a homogeneous orthotropic fiber reinforced composite with coinciding fiber directions. 
Under the assumptions of perfect bonding between fiber and matrix material, small deformations, and a uniform stress distribution in the fiber, the displacement and stress continuity conditions at the fiber matrix interface are evaluated using expressions from \cite{Lekhnitskii:1981:1} to describe the linearly elastic effect of the elliptic inclusion in the matrix material. Residual strains in the fiber are neglected and the fiber is assumed to have infinite length.
This leads to a 1:1 relation between $(\varepsilon(\uf))_1$ and $(\varepsilon(\unf))_1$, for $\varepsilon(\uf)$ and $\varepsilon(\unf)$ in Voigt notation, \ie,
\begin{equation*}
	\varepsilon(u)
	= 
	\begin{pmatrix}
		\varepsilon(u)_{11} &
		\varepsilon(u)_{22} &
		\varepsilon(u)_{33} &
		2 \varepsilon(u)_{23} &
		2 \varepsilon(u)_{13} &
		2 \varepsilon(u)_{12} 
	\end{pmatrix}^\transp
	.
\end{equation*}
We assume here that the first axis of the coordinate system is aligned with the fiber direction.
We additionally assume constant temperatures, \ie $\Delta T=0$. 
The detailed derivation in \cite{KollarVanSteenkiste:1998:1} then results in a sparse strain transfer matrix $T$ with the non-zero entries expressed by the relations
\begin{subequations}\label{eq:analytical_stp}
	\begin{align}
		\label{eq:analytical_stp_1}
		(\varepsilon(\uf))_1 
		&
		= 
		(\varepsilon(\unf))_1,
		\\
		\label{eq:analytical_stp_2}
		\begin{pmatrix}
			(\varepsilon(\uf))_2
			\\
			(\varepsilon(\uf))_3
			\\
			(\varepsilon(\uf))_4
			\\
			\Theta(\uf)
		\end{pmatrix}
		&
		= 
		\paren[big](){O-UL^{-1}N}^{-1} 
		\paren[auto][.{%
			U L^{-1} W 
			\begin{pmatrix}
				\paren[big](){C_{21}^f - C_{21}} (\varepsilon(\unf))_1\\ \paren[big](){C_{31}^f - C_{31}} (\varepsilon(\unf))_1\\0
			\end{pmatrix}
		}
		\breakeqn{1}
		+ 
		\paren[auto].]{%
			\paren[big](){K-U L^{-1} W Q}
			\begin{pmatrix}
				(\varepsilon(\unf))_2
				\\
				(\varepsilon(\unf))_3
				\\
				(\varepsilon(\unf))_4
			\end{pmatrix} 
		}
		,
		\\
		\label{eq:analytical_stp_3}
		(\varepsilon(\uf))_5 
		&
		=
		\frac{b+\sqrt{\frac{C_{55}}{C_{66}}} a}{b+\frac{C^f_{55}}{\sqrt{C_{66}C_{55}}}a} (\varepsilon(\unf))_5,
		\\
		\label{eq:analytical_stp_4}
		(\varepsilon(\uf))_6 
		&
		=
		\frac{a+\sqrt{\frac{C_{66}}{C_{55}}} b}{a+\frac{C^f_{66}}{\sqrt{C_{66}C_{55}}}b} (\varepsilon(\unf))_6
		.
	\end{align}
\end{subequations}
Here $\Theta(\uf)$ is the angular displacement of the sensor, which we ignore here.
Moreover, we denote the entries of the stiffness tensor $\C$ in Voigt notation by $C_{ij}$, $1 \leq i,j \leq 6$ and, similarly, the entries of the compliance tensor $\S=\C^{-1}$ in Voigt notation by $S_{ij}$, $1 \leq i,j \leq 6$. 
The superscript $\cdot\,^f$ denotes material parameters of the fiber and no superscript denotes matrix material parameters. 
Furthermore, $a = b = R$ denotes the lengths of the semi-axes of the fiber's cross-section.
The matrices in \eqref{eq:analytical_stp_2} are given by
\begin{equation*}
	K 
	= 
	\begin{pmatrix}
		a & 0 & 0\\
		0 & 0 & \frac{b}{2}\\
		0 & 0 & \frac{a}{2}\\
		0 & b & 0
	\end{pmatrix}
	, 
	\quad 
	H 
	= 
	\begin{pmatrix}
		0\\ b\\ -a\\ 0
	\end{pmatrix}
	,
	\quad 
	W 
	= 
	\begin{pmatrix}
		1 & 0 & 0\\
		0 & 1 & 0\\
		0 & 0 & 1\\
		0 & 0 & 1
	\end{pmatrix}
	, 
\end{equation*}
\begin{equation*}
	Q 
	= 
	\begin{pmatrix}
		C_{22} & C_{23} & 0\\
		C_{32} & C_{33} & 0\\
		0 & 0 & C_{44}
	\end{pmatrix}
	,
	\quad
	U 
	= 
	\begin{pmatrix}
		\delta_1 & -\delta_2 &  \delta_1 &  \delta_2 \\
		\delta_2 &  \delta_1 & -\delta_2 &  \delta_1 \\
		\delta_3 & -\delta_4 & -\delta_3 & -\delta_4 \\
		\delta_4 &  \delta_3 &  \delta_4 & -\delta_3
	\end{pmatrix},
\end{equation*}
with
\begin{equation*}
	\begin{aligned}
		\delta_1 
		&
		=
		2 \beta_{23} + 2 \beta_{22} (\mu_R^2 - \mu_I^2)
		,
		\quad
		\delta_2
		=
		4 \beta_{22} \mu_R \mu_I
		,
		\quad
		\delta_3
		=
		2 \mu_R \left( \beta_{23} + \frac{\beta_{33}}{\mu_R^2+\mu_I^2} \right)
		,
		\quad
		\breakeqn{2}
		\delta_4
		&
		=
		2 \mu_I \left( \beta_{23} - \frac{\beta_{33}}{\mu_R^2+\mu_I^2} \right)
	\end{aligned}
\end{equation*}
and
\begin{equation*}
	L = \begin{pmatrix}
		\frac{2 \mu_I}{b} & \frac{2 \mu_R}{b} & \frac{2 \mu_I}{b} & -\frac{2 \mu_R}{b}\\
		\frac{2}{a} & 0 & \frac{2}{a} & 0\\
		0 & -\frac{2}{b} & 0 & -\frac{2}{b}\\
		-\frac{2 \mu_R}{a} & \frac{2 \mu_I}{a} & \frac{2 \mu_R}{a} & \frac{2 \mu_I}{a}
	\end{pmatrix}, \quad    
	O=\begin{pmatrix}
		K & H
	\end{pmatrix}, \quad
	N=\begin{pmatrix}
		W Q_f & \boldsymbol{0}
	\end{pmatrix},
\end{equation*}
where $\bnull$ denotes the zero vector $\bnull \in \mathbb{R}^{4 \times 1}$.
The matrix $Q_f$ is similar to $Q$ but with fiber material parameters $C^f_{ij}$ instead of $C_{ij}$. 
Finally, by \cite{Lekhnitskii:1981:1}, the parameters $\beta_{ij}, \mu_R$ and $\mu_I$ are related to the entries of $\S$ via 
\begin{equation*}
	\begin{aligned}
		\beta_{ij}
		&
		=
		S_{ij} - \frac{S_{i1}S_{j1}}{S_{11}}
		, 
		\quad 
		\mu_R
		=
		\sqrt{\sqrt{\frac{\beta_{33}}{4 \beta_{22}}} - \frac{2 \beta_{23} + \beta_{44}}{4 \beta_{22}}}
		, 
		\quad 
		\breakeqn{2}
		\mu_I
		&
		=
		\sqrt{\sqrt{\frac{\beta_{33}}{4 \beta_{22}}} + \frac{2 \beta_{23} + \beta_{44}}{4 \beta_{22}}}
		.
	\end{aligned}
\end{equation*}

\section{Extension of the Strain Transfer Principle}
\label{sec:strain_transfer_principle_extensions}

We now consider the case $E_f > E_m$, which is relevant for fiber Bragg grating applications.
A typical case is that of a glass fiber embedded in either an isotropic or fiber reinforced plastic.
For the sake of simplicity, we ignore here the soft protective coating of the fiber.
We propose the following extension of the STP (in tensor notation)
\begin{equation}
	\label{eq:extended_strain_transfer_principle}
	\varepsilon(\uf) \approx (\I - \projgamma \otimes \projgamma) : \left( \Tnf : \varepsilon(\unf) \right) + \projgamma \varepsilon_\gamma(\uf) \quad \text{on } \gamma.
\end{equation}
\Cref{eq:extended_strain_transfer_principle} lets us recover the strain tensor $\varepsilon(\uf)$ using the strain field $\varepsilon(\unf)$ computed from the solution~$\unf$ of \eqref{eq:min-problem_without_fiber} (without a fiber), and using the strain in fiber direction $\varepsilon_\gamma(\uf)$ computed from the solution~$\uf$ of \eqref{eq:min-problem} (taking the stiffening due to the fiber into account).
Note that since we consider a one-dimensional fiber, the full $\varepsilon(\uf)$ tensor is not directly computable from $\uf$ in a meaningful way, but only the component~$\varepsilon_\gamma(\uf)$ in fiber direction is available. 
We remark that \eqref{eq:extended_strain_transfer_principle} is equivalent to the original STP, but with the strain in fiber direction corrected.

\section{Variational Formulation and Numerical Discretization}\label{sec:variational_form}

Problem \eqref{eq:min-problem} will be solved using the method of finite elements.
The displacement~$u \in \Uad$  minimizing the energy in \eqref{eq:min-problem} is characterized by the variational formulation 
\begin{equation}
	\label{eq:variational_form}
	\int_{\Omega} \varepsilon(u) \dprod \C \dprod \varepsilon(v) \dx + \int_{\gamma} E A \, \varepsilon_\gamma(u) \, \varepsilon_\gamma(v) \ds 
	= 
	\int_{\Omega} f \cdot v \dx \quad \text{for all } v \in V,
\end{equation}
where $\Uad = \setDef{u \in H^1(\Omega)}{u = u_0 \text{ on } \Gamma_D}$ and $\Gamma_D$ is the boundary of $\Omega$ with imposed Dirichlet boundary conditions~$u = u_0$ for the displacement. 
The corresponding test space $V$ is given by $V = \setDef{u \in H^1(\Omega)}{u = 0 \text{ on } \Gamma_D}$.
For the finite element discetization, the domain $\Omega \subset \R^3$ is approximated by a tetrahedral mesh consisting of a set of tetrahedrons $\mathcal{T}(\Omega)$ such that the fiber~$\gamma$ is approximated by a set of edges $\mathcal{E}(\gamma)$ of the mesh.
The set of admissible displacements $\Uad$ is approximated by functions $\Uadh$ which are piecewise linear on each tetrahedron, globally continuous and satisfy the Dirichlet boundary conditions. 
Similarly, the set of test functions is approximated by functions $V_h$ which are piecewise linear on each tetrahedron, globally continuous and are zero on the boundary $\Gamma_D$. 
Due to the linearity, $\varepsilon(u)$ and $\varepsilon(v)$ are constant on each tetrahedron.
In this discrete setting, the variational form \eqref{eq:variational_form} becomes
\begin{align}
	\sum_{t \in \mathcal{T}(\Omega)}  \int_{t} \varepsilon(u) \dprod \C \dprod \varepsilon(v) \dx + \sum_{e \in \mathcal{E}(\gamma)} \int_{e} E A  \varepsilon_\gamma(u) \varepsilon_\gamma(v) \ds 
	= 
	\sum_{t \in \mathcal{T}(\Omega)}  \int_{t} f \cdot v \dx 
	\quad 
	\breakeqn{3}
	\text{for all } v \in V_h
	.
	\label{eq:discrete_variational_form}
\end{align}
For an edge $e \in \mathcal{E}(\gamma)$ and its incident vertices $p_1$ and $p_2$, the strain in fiber direction can be calculated as
\begin{equation}
	\varepsilon_\gamma(u) = \frac{(u(p_2) - u(p_1)) \cdot (p_2 - p_1)}{|p_2 - p_1|^2},
\end{equation}
where $u(p_1)$ and $u(p_2)$ denote the nodal displacements in $p_1$ and $p_2$.

\section{Numerical Demonstration}
\label{sec:experiments}

We test the proposed extension of the STP on an example of intermediate complexity.
The geometry is a $\SI{100}{\milli\meter} \times \SI{100}{\milli\meter} \times \SI{10}{\milli\meter}$ plate with a $\SI{32}{\milli\meter}$ diameter bore located at $(x,y) = (\SI{60}{\milli\meter}, \SI{40}{\milli\meter})$, as shown in \cref{fig:domain}. 
As matrix material we use a fiber reinforced plastic, for which we compute the effective stiffness tensor by homogenization, described below in \cref{sec:homogenization}.
The sensor fiber is made of glass.
All FE computations were performed using \dolfin/\fenics 2019.1; see \cite{AlnaesBlechtaHakeJohanssonKehletLoggRichardsonRingRognesWells:2015:1,LoggMardalWells:2012:1}.
The strain transfer described in \eqref{eq:analytical_stp} is also computed numerically. 
Further details are given in the following sections.

\subsection{Domain and Mesh Generation}
\label{sec:meshing}

\begin{figure}[htbp]
	\begin{subfigure}[c]{0.5\textwidth}
		\includegraphics[height=0.9\linewidth,trim=100 95 40 105,clip]{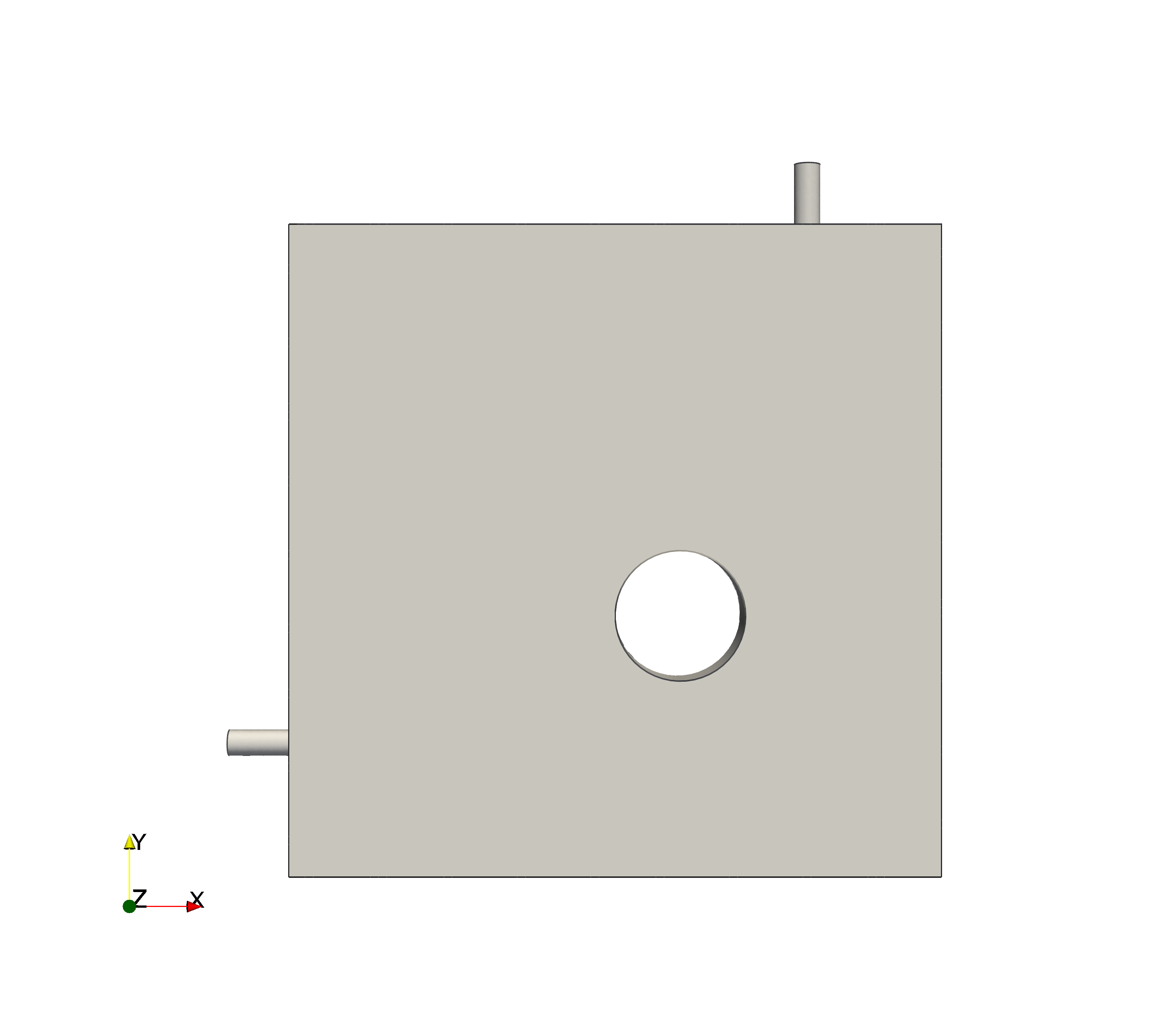}%
		\subcaption{domain $\Omega$}
		\label{fig:domain}
	\end{subfigure}
	\begin{subfigure}[c]{0.5\textwidth}
		\includegraphics[height=0.9\linewidth,trim=50 30 10 30,clip]{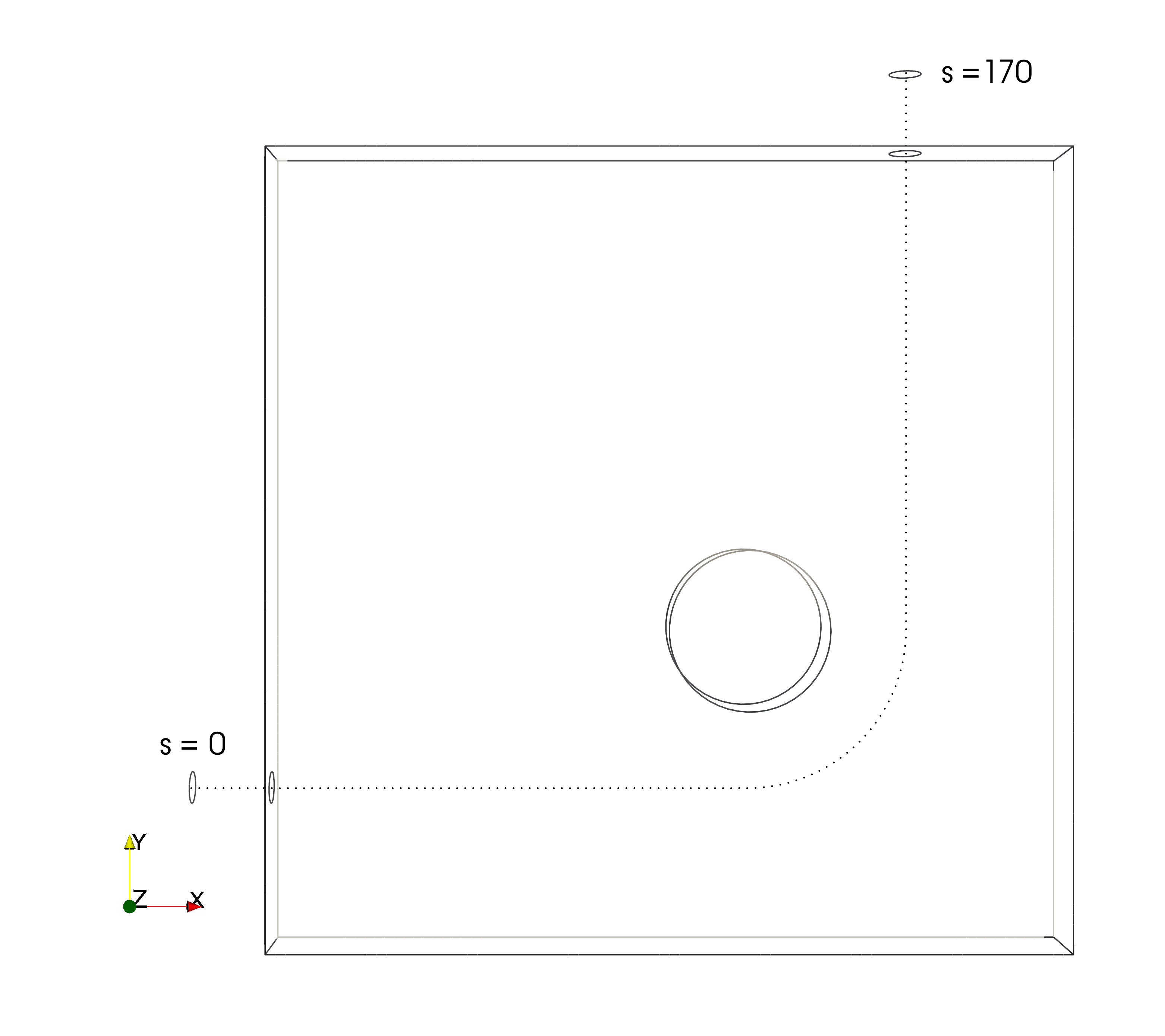}%
		\subcaption{path of the fiber center}
		\label{fig:fiber_path}
	\end{subfigure}
	\begin{subfigure}[c]{0.5\textwidth}
		\includegraphics[height=0.9\linewidth,trim=50 30 10 30,clip]{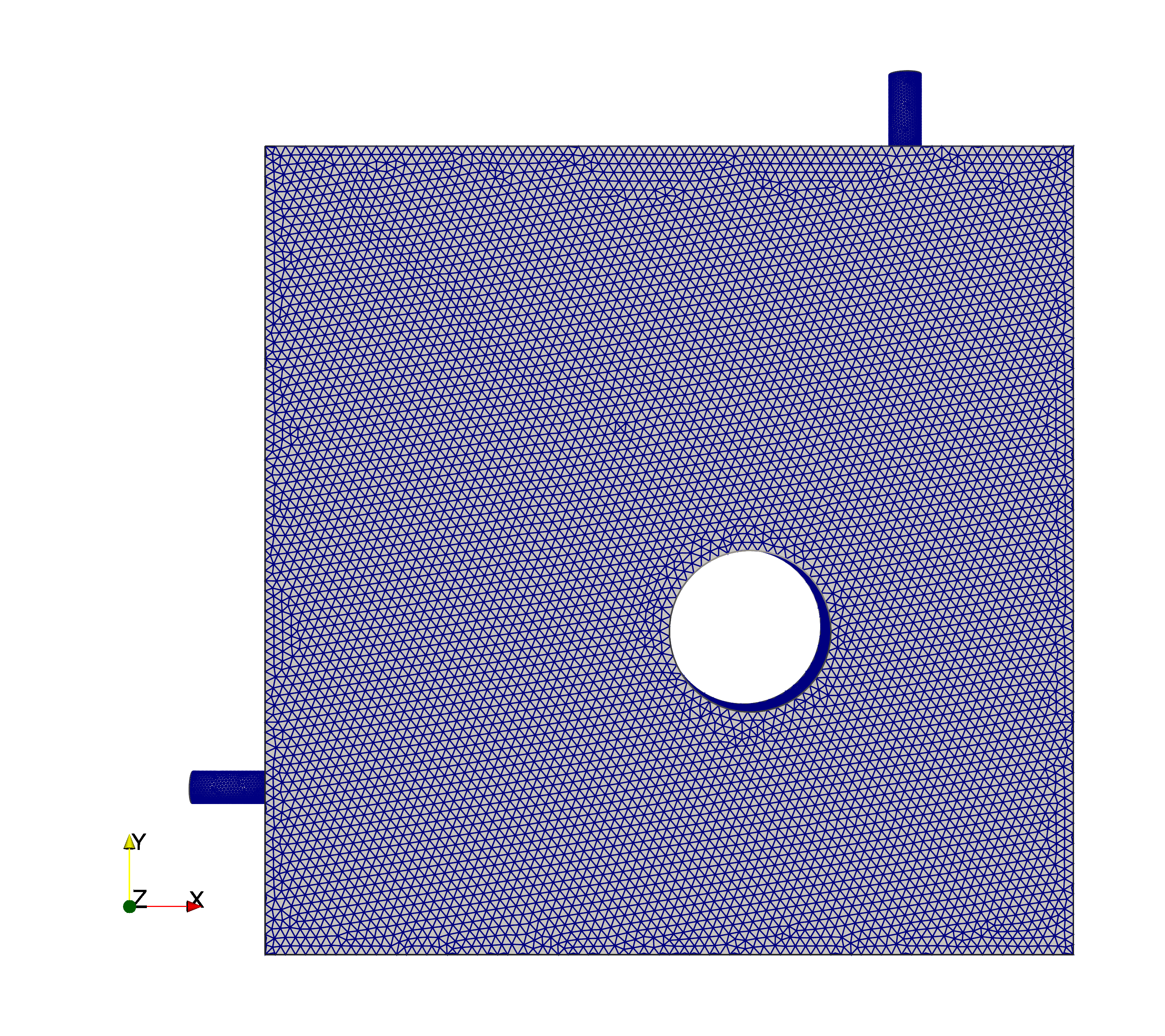}%
		\subcaption{mesh on the boundary}
		\label{fig:boundary_mesh}
	\end{subfigure}
	\begin{subfigure}[c]{0.5\textwidth}
		\includegraphics[height=0.9\linewidth,trim=50 30 10 30,clip]{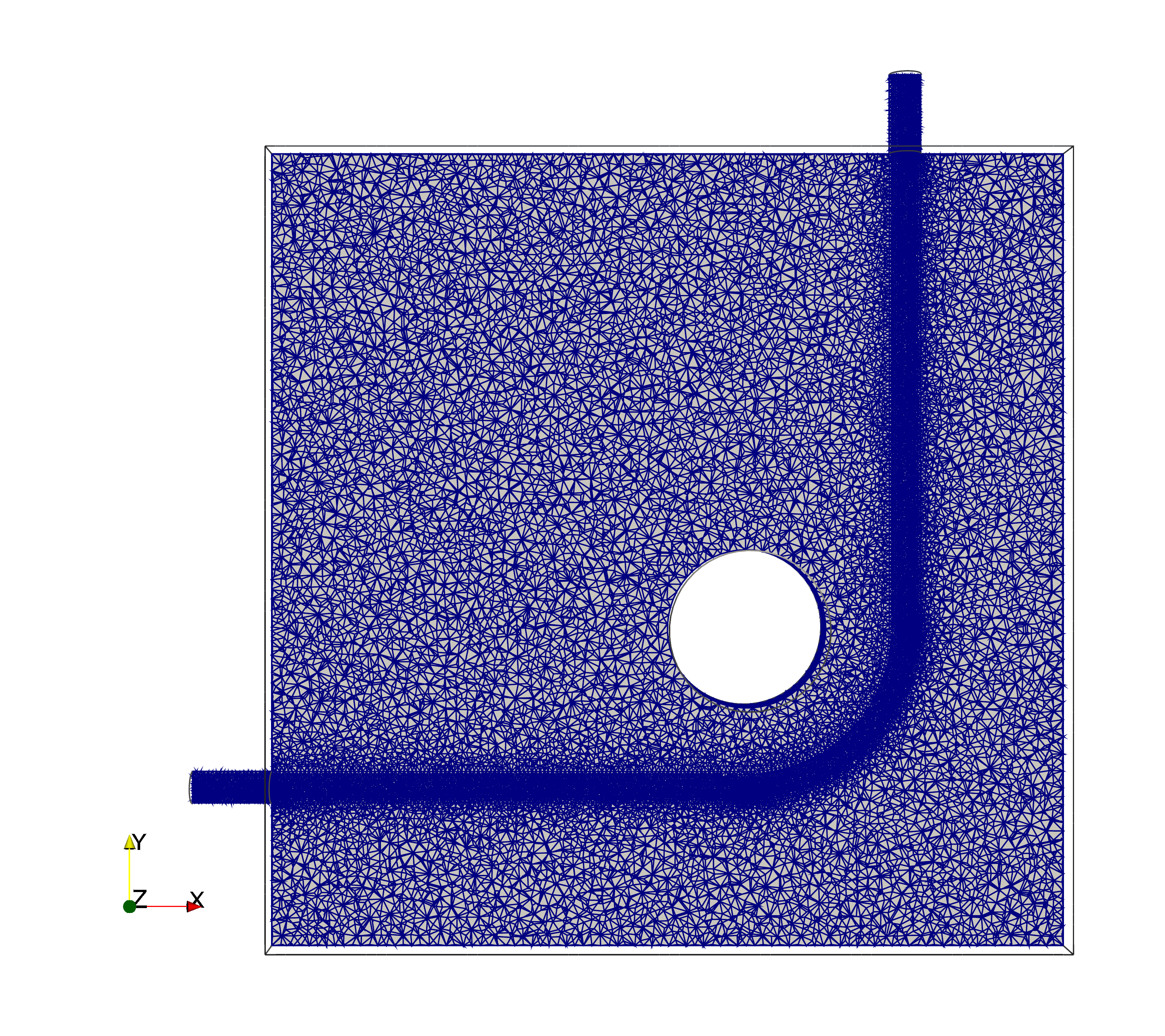}%
		\subcaption{interior mesh}
		\label{fig:inner_mesh}
	\end{subfigure}
	\caption{Plate with hole \subref{fig:domain}, path of the sensor fiber \subref{fig:fiber_path} and generated mesh \subref{fig:boundary_mesh}, \subref{fig:inner_mesh}.}
	\label{fig:mesh}
\end{figure}

The geometry shown in \cref{fig:domain} was created using \freecad\footnote{\url{https://www.freecadweb.org}}.
Notice that this geometry also contains a fiber of diameter $\SI{4}{\milli\meter}$ and length $s= \SI{170}{\milli\meter}$ along the path shown in \cref{fig:fiber_path}, which begins at $(x,y,z) = (\SI{-10}{\milli\meter}, \SI{20}{\milli\meter}, \SI{5}{\milli\meter})$ and ends at $(x,y,z) = (\SI{80}{\milli\meter}, \SI{110}{\milli\meter}, \SI{5}{\milli\meter})$. 
The sole purpose of resolving the relatively thick fiber in the mesh is to create a reference finite element solution to compare to the results obtained by the STP.
The fiber cross-section itself was split into 4~quadrants such that after meshing there will be an edge following the center of the fiber. For meshing, the geometry was exported from \freecad as a STEP file. 
The STEP file was then loaded into \gmsh using the \opencascade plugin. 
The boundaries of the different regions were associated using the \eqq{Coherence} function in \gmsh and different subdomains, boundaries and the fiber center line were labeled.
The characteristic length used for meshing was calculated from curvature, such that the mesh is more refined near the fiber.
The resulting mesh at the boundary and inside of the domain can be seen in \cref{fig:boundary_mesh} and \cref{fig:inner_mesh}, respectively. 
The mesh contains $\num{161476}$ nodes and $\num{963272}$ tetrahedral elements. 
The \gmsh mesh including subdomains, boundaries and paths was then loaded and converted to the XDMF format using \meshio\footnote{\url{https://github.com/nschloe/meshio}}.

As was mentioned above, the mesh with the three-dimensional fiber resolved is used for the purpose of computing reference solutions. 
However, the same mesh was also used when computing the strains for the embedded one-dimensional fiber using the STP. 
In this case, the material inside of the meshed fiber was set equal to the matrix material.
We followed this procedure to avoid the influence of different meshes on the solutions.
In practice, there would be no need to refine the mesh close to the fiber, nor to resolve the fiber in the mesh.

\subsection{Homogenization of Matrix Material}\label{sec:homogenization}

For the demonstration, we used glass fiber reinforced polypropylene as matrix material. 
For glass we use a Young's modulus of $E_g = \SI{73}{\giga\pascal}$ and a Poisson ratio of $\nu_g = 0.18$. 
For polypropylene we use $E_{pp} = \SI{1.665}{\giga\pascal}$ and $\nu_{pp} = 0.36$. 
The fibers are assumed to be parallel and of infinite length in $z$-direction with a diameter which results in a fiber volume fraction of $50\%$.
The effective stiffness matrix was computed using \fibergen from \cite{Ospald:2019:1}, which employs a Lippmann-Schwinger approach (see \cite{MoulinecSuquet:1998:1}) with a conjugate gradient solver on a staggered grid described in \cite{KabelBoehlkeSchneider:2014:1,SchneiderOspaldKabel:2016:1}. 
Using a laminate mixing rule at the interfaces, see \cite{SchneiderOspaldKabel:2016:2}, only a resolution of $55 \times 32 \times 1$ voxels is required for the representative volume element to achieve a sufficient accuracy.
The method also allows the computation of effective material properties for other fiber distributions, \eg, for injection molded parts.

Using these settings, the homogenized matrix material stiffness tensor in Voigt notation reads
\begin{equation}
	\C = \begin{pmatrix}
		6.34 & 3.03 & 2.43 & 0 & 0 & 0 \\
		3.03 & 6.34 & 2.43 & 0 & 0 & 0 \\
		2.43 & 2.43 & 38.61& 0 & 0 & 0 \\
		0 & 0 & 0 & 1.75 & 0 & 0 \\
		0 & 0 & 0 & 0 & 1.75 & 0 \\
		0 & 0 & 0 & 0 & 0 & 1.65
	\end{pmatrix} \SI{}{\giga\pascal},
\end{equation}
which enters the computation of the strain transfer matrix in the following section.
For our FE calculations, we assume the reinforcement fibers to be oriented in the $x$-direction of our plate. In this instance, one has to swap the $x$ with the $z$-axis of $\C$ to obtain the correct material law.

\subsection{Computation of Strain Transfer Matrices}\label{sec:computation_ST_matrices}

For the embedded fiber, we assume the same properties of glass as above, \ie $E_f = \SI{73}{\giga\pascal}$ and a Poisson ratio of $\nu_f = 0.18$. 
The strain transfer matrix between strain tensors in Voigt notation representing the respective strain transfer tensor~$T$ in \eqref{eq:strain_transfer_principle} can then be computed analytically as demonstrated in \Cref{sec:analytical_stp}.
Alternatively, we can evaluate it numerically.
We chose the latter approach using \fibergen \cite{Ospald:2019:1} in a similar fashion as for the homogenization.
As representative volume element for the latter, we chose a $1 \times 1$ box with a disc of diameter $0.05$ placed at the center representing the fiber, while the remaining domain represents the matrix material.
Due to use of periodic boundary conditions for the displacements, the diameter of the disc has to be sufficiently small and the resolution sufficiently large (in our case $512 \times 512$ voxels).
For the identification of the strain transfer matrix, six linearly independent load cases with prescribed strain $E^{(i)}$ are required. The prescribed strain represents the far field or matrix strain at the fiber position.
The computed strain field $\varepsilon^{(i)}$ for prescribed strain $E^{(i)}$
is evaluated at the center of the domain to obtain the strain inside the fiber $\varepsilon^{(i)}_f = \varepsilon^{(i)}(0.5, 0.5)$.
The transfer matrix $T$  is then given by the relation
\begin{equation*}
	\begin{pmatrix}
		\varepsilon^{(1)}_f | \dots | \varepsilon^{(6)}_f
	\end{pmatrix} = T \begin{pmatrix}
		\varepsilon^{(1)}_m | \dots | \varepsilon^{(6)}_m
	\end{pmatrix},
\end{equation*}
where the strains in each column are given in Mandel notation (the notation is only important to interpret the numerical values below).
In our instance, the fiber is oriented in $z$-direction, \ie parallel to the reinforcement fibers of the matrix material. 
In order to compute the strain transfer matrix for instances where the sensor fiber has an angle $\alpha$ to the reinforcement fibers we keep the sensor fiber oriented in $z$-direction but rotate the matrix material around the $x$-axis.
The rotated $\C^R$ (in full tensor notation) is then given by
\begin{equation*}
	\C^R_{ijk\ell} = R_{im} R_{jn} R_{kp} R_{\ell q} \C_{mnpq},
\end{equation*}
where the sum is carried out over all free indices (using Einstein summation) and $R=R_x$ represents the rotation matrix for rotation by the angle $\alpha$ around the $x$-axis
\begin{equation*}
	R_x = \begin{pmatrix}
		1 & 0 & 0 \\
		0 & \cos \alpha & -\sin \alpha \\
		0 & \sin \alpha & \cos \alpha
	\end{pmatrix}.
\end{equation*}
Similarly as for $\C$ of our plate, one has to swap the $x$- with the $z$-axis of $T$ to obtain the strain transfer matrix in the correct coordinate system for our example. Furthermore, if the fiber orientation (in the $x$-$y$-plane) has an angle $\beta \neq 0$ to the $x$-axis one has to rotate $T$ around the $z$-axis by angle $\beta$, \ie
\begin{equation*}
	T^R_{ijk\ell} = R_{im} R_{jn} R_{kp} R_{\ell q} T_{mnpq},
\end{equation*}
where again the sum is carried out over all free indices (Einstein summation) and $R=R_z$ represents the rotation matrix for a rotation by the angle $\beta$ around the $z$-axis
\begin{equation*}
	R_z = \begin{pmatrix}
		\cos \beta & -\sin \beta & 0 \\
		\sin \beta & \cos \beta & 0 \\
		0 & 0 & 1
	\end{pmatrix}.
\end{equation*}
Note that $T$ also has to be blown up to a full $4$-tensor and then converted back to a matrix in Mandel notation.

In the first horizontal section of the path of the sensor fiber in our example (see \cref{fig:fiber_path}) we have $\beta = 0$ and the computed (and properly rotated) strain transfer matrix is given by
\begin{equation*}
	T = \begin{pmatrix} 
		1.00 & 0 & 0 & 0 & 0 & 0 \\
		-0.15 & 0.10 & 0.02 & 0 & 0 & 0 \\
		-0.15 & 0.02 & 0.10 & 0 & 0 & 0 \\
		0 & 0 & 0 & 0.08 & 0 & 0 \\
		0 & 0 & 0 & 0 & 0.11 & 0 \\
		0 & 0 & 0 & 0 & 0 & 0.11
	\end{pmatrix}.
\end{equation*}
In the middle of the arc section at an angle of $\beta = \frac{\pi}{4}$, the strain transfer matrix is given by
\begin{equation*}
	T = \begin{pmatrix} 
		0.67 & 0.16 & 0.01 & 0 & 0 & 0.25 \\
		0.03 & 0.30 & 0.01 & 0 & 0 & 0.27 \\
		-0.11 & -0.06 & 0.11 & 0 & 0 & -0.12 \\
		0 & 0 & 0 & 0.10 & 0.03 & 0 \\
		0 & 0 & 0 & 0.03 & 0.13 & 0 \\
		0.21 & 0.38 & -0.01 & 0 & 0 & 0.63
	\end{pmatrix}
\end{equation*}
and finally in the vertical section ($\beta = \frac{\pi}{2}$) we have
\begin{equation*}
	T = \begin{pmatrix} 
		0.55 & -0.14 & 0.02 & 0 & 0 & 0 \\
		0 & 1.00 & 0 & 0 & 0 & 0 \\
		-0.08 & -0.14 & 0.11 & 0 & 0 & 0 \\
		0 & 0 & 0 & 0.10 & 0 & 0 \\
		0 & 0 & 0 & 0 & 0.11 & 0 \\
		0 & 0 & 0 & 0 & 0 & 0.11
	\end{pmatrix}.
\end{equation*}
Note that in all sections the strain in fiber direction is always transferred verbatim from the matrix material, as recognized, \eg, from the unit diagonal entry in the first and the last fiber sections.

\subsection{Solution of Linear Elasticity Problems}

The discretization and solution of \eqref{eq:variational_form}, its counterpart coming from \eqref{eq:matrix_energy}, and the reference solution are performed using \dolfin/\fenics~2019.1; see \cite{AlnaesBlechtaHakeJohanssonKehletLoggRichardsonRingRognesWells:2015:1,LoggMardalWells:2012:1}.
Meshes, subdomains and the fiber path are loaded from the XDMF files generated as described in \cref{sec:meshing}. 
The plate is clamped on the lower boundary ($y=\SI{0}{\milli\metre}$, see \cref{fig:domain}) and a fixed displacement of $(1,0,0)^\transp$~\SI{}{\milli\metre} is enforced on the upper ($y=\SI{100}{\milli\metre}$) boundary.
The volume force $f$ is set to zero.
For the solution of the arising linear systems we use the conjugate gradient method together with an AMG preconditioner from the \dolfin \petsc backend.
Three solutions are obtained: the displacement field of the reference solution $u_r$ (with three-dimensionally resolved fiber); the solution $\unf$ with the fiber neglected by setting the material inside the fiber subdomain to the matrix material; and the solution $\uf$ from the superimposed one-dimensional fiber model \eqref{eq:variational_form}, where the material inside the fiber subdomain is also set to the matrix material but the fiber's stiffness enters through the stretching energy term.

\subsection{Evaluation}

\newcommand{\plotstrains}[6]{
	\tikzsetnextfilename{#3_strain_#1}
	\hfill
	\begin{tikzpicture}[scale=1]
		\begin{axis}[
			xlabel=$s\,(\si{\milli\meter})$,
			ylabel=$\varepsilon_{#1}$,
			xmin=0.0,
			xmax=150.0,
			ymin=#4,
			ymax=#5,
			legend columns=#6,
			grid=major,
			grid style={dashed,gray!30},
			legend pos=north east,
			legend cell align=left,
			width=0.95\linewidth
		]
			\addplot[color=black, mark=star, mark options={solid, black, fill=white, fill opacity=0}] table[x=s,y expr=\thisrow{er#1}] {#2};
			\addlegendentry{$\varepsilon^{r}_{#1}$}
			\addplot[color=green, mark=square, mark options={solid, green, fill=white, fill opacity=0}] table[x=s,y expr=\thisrow{enf#1}] {#2};
			\addlegendentry{$\varepsilon^{nf}_{#1}$}
			\addplot[color=red, mark=triangle, mark options={solid, red, fill=white, fill opacity=0}] table[x=s,y expr=\thisrow{eTnf#1}] {#2};
			\addlegendentry{$T\varepsilon^{nf}_{#1}$}
			\addplot[color=blue, mark=diamond, mark options={solid, blue, fill=white, fill opacity=0}] table[x=s,y expr=\thisrow{ef#1}] {#2};
			\addlegendentry{$\varepsilon^{f}_{#1}$}
		\end{axis}
	\end{tikzpicture}
}

\pgfplotstableread{results/plate_shear_top_0.csv}{\PlateResultsZero}
\pgfplotstableread{results/plate_shear_top_50.csv}{\PlateResultsFifty}

\begin{figure}[htbp]
	\begin{subfigure}[c]{0.5\textwidth}
		\plotstrains{11}{\PlateResultsZero}{plate_0}{-2e-3}{3e-3}{1}%
		\subcaption{}
		\label{fig:results_11_0}
	\end{subfigure}
	\begin{subfigure}[c]{0.5\textwidth}
		\plotstrains{11}{\PlateResultsFifty}{plate_50}{-2e-3}{3e-3}{1}%
		\subcaption{}
		\label{fig:results_11_50}
	\end{subfigure}
	\newline
	\begin{subfigure}[c]{0.5\textwidth}
		\plotstrains{22}{\PlateResultsZero}{plate_0}{-2.5e-3}{6e-3}{2}%
		\subcaption{}
		\label{fig:results_22_0}
	\end{subfigure}
	\begin{subfigure}[c]{0.5\textwidth}
		\plotstrains{22}{\PlateResultsFifty}{plate_50}{-2.5e-3}{6e-3}{2}%
		\subcaption{}
		\label{fig:results_22_50}
	\end{subfigure}
	\newline
	\begin{subfigure}[c]{0.5\textwidth}
		\plotstrains{33}{\PlateResultsZero}{plate_0}{-3e-3}{3e-3}{2}%
		\subcaption{}
		\label{fig:results_33_0}
	\end{subfigure}
	\begin{subfigure}[c]{0.5\textwidth}
		\plotstrains{33}{\PlateResultsFifty}{plate_50}{-3e-3}{3e-3}{2}%
		\subcaption{}
		\label{fig:results_33_50}
	\end{subfigure}
	\caption{Main strain components of the solutions along the path of the sensor fiber. The plots on the left show the results for a matrix fiber volume fraction of $0\%$. The plots on the right are for a matrix fiber volume fraction of $50\%$.}
\label{fig:results_1}
\end{figure}

\begin{figure}[htbp]
	\begin{subfigure}[c]{0.5\textwidth}
		\plotstrains{23}{\PlateResultsZero}{plate_0}{-2e-4}{3.5e-4}{2}%
		\subcaption{}
		\label{fig:results_23_0}
	\end{subfigure}
	\begin{subfigure}[c]{0.5\textwidth}
		\plotstrains{23}{\PlateResultsFifty}{plate_50}{-2e-4}{3.5e-4}{2}%
		\subcaption{}
		\label{fig:results_23_50}
	\end{subfigure}
	\newline
	\begin{subfigure}[c]{0.5\textwidth}
		\plotstrains{13}{\PlateResultsZero}{plate_0}{-2e-4}{3.5e-4}{2}%
		\subcaption{}
		\label{fig:results_13_0}
	\end{subfigure}
	\begin{subfigure}[c]{0.5\textwidth}
		\plotstrains{13}{\PlateResultsFifty}{plate_50}{-2e-4}{3.5e-4}{2}%
		\subcaption{}
		\label{fig:results_13_50}
	\end{subfigure}
	\newline
	\begin{subfigure}[c]{0.5\textwidth}
		\plotstrains{12}{\PlateResultsZero}{plate_0}{-0.5e-3}{8.5e-3}{2}%
		\subcaption{}
		\label{fig:results_12_0}
	\end{subfigure}
	\begin{subfigure}[c]{0.5\textwidth}
		\plotstrains{12}{\PlateResultsFifty}{plate_50}{-0.5e-3}{8.5e-3}{2}%
		\subcaption{}
		\label{fig:results_12_50}
	\end{subfigure}
	\caption{Mixed strain components of the solutions along the path of the sensor fiber. The plots on the left show the results for a matrix fiber volume fraction of $0\%$. The plots on the right are for a matrix fiber volume fraction of $50\%$.}
	\label{fig:results_2}
\end{figure}

In \cref{fig:results_1,fig:results_2} we plot all components of the corresponding strain tensors along the path of the fiber. 
The strain for the reference solution~$u_r$ is denoted by $\varepsilon^r$, the strain for $\unf$ is given by $\varepsilon^{nf}$ and the strain $\varepsilon^f$ for $\uf$ is obtained from our extended STP given in \eqref{eq:extended_strain_transfer_principle}.
For comparison, the strain obtained by the original STP \eqref{eq:strain_transfer_principle} is denoted by $T\varepsilon^{nf}$.
In addition to a matrix fiber volume fraction of $50\%$ (right plots of \cref{fig:results_1,fig:results_2}), we also performed the same simulations with a fiber volume fraction of $0\%$ (left plots of \cref{fig:results_1,fig:results_2}) representing a very soft and isotropic matrix material.

Ideally, the STP solutions $T\varepsilon^{nf}$ (red triangles) and $\varepsilon^f$ (blue diamonds) should be identical to the reference solution $\varepsilon^r$ (black stars).
For the $11$- and $22$-components we can see that our extended version of the STP agrees very well with the reference solution, whereas the original STP has major deviations for the $11$-component until after the bend of the sensor fiber at around $s=\SI{90}{\milli\meter}$ as well as for the $22$-component beginning with the bend at around $s=\SI{60}{\milli\meter}$.
This observation is independent of the matrix fiber volume fraction.
For the $33$-component the extended STP and original STP are identical, since the sensor fiber direction is always orthogonal to the $z$-direction. 
Also they both disagree with the reference solution by a significant amount for the $0\%$ matrix fiber volume fraction case.
The $23$- and $13$-components are one order of magnitude smaller than the other components and in theory they should be zero in view of the symmetry of the problem in $z$-direction.
Finally, for the $12$-component the original and extended STP agree everywhere except for the bend around the hole. 
Here again the extended STP outperforms the original STP in the case of a $0\%$ matrix fiber volume fraction. 
The $50\%$ case is indecisive.

\section{Conclusion and Outlook}
\label{sec:conclusion_outlook}

In this paper, we proposed an improvement of the original strain transfer principle, which recovers the strain components inside a fiber embedded in a matrix material from simulations which do not require the fiber geometry to be resolved.
The matrix material itself can be isotropic or fiber reinforced.
The proposed modification to the classical STP consists of an additional term in the elastic energy of the total part, which takes into account the additional stretching energy using a simple one-dimensional fiber model.
This modification is particularly relevant for fiber materials which are stiffer than the matrix, as it is often the case for fiber optical strain sensors made of glass.
Our evaluation shows that the extended STP improves the classical STP in regions in which the presence of the sensor fiber restricts the displacement of the surrounding material in fiber direction. 

It is a limitation that the one-dimensional fiber model \eqref{eq:min-problem} does not incorporate lateral strains, which results in moderate deviations in the $33$-component. 
More deviations are expected due to the neglection of bending and twisting terms.
While these have only a minor contribution to the overall energy of the fiber in the chosen example, they may become more relevant in other setups.

The inclusions of bending and twisting energies as well as the consideration of lateral strains are left to future research.
We expect that these terms will be quite challenging to model, discretize and implement.
Furthermore, as glass fiber sensors are usually coated with a protective layer made of a soft material, a model for coated fibers should be considered. 
For this case there already exists an analytical STP proposed in the literature; see for instance \cite{VanSteenkisteKollar:1998:1}. 
Additionally, the integration of ideas from the shear-lag theory mentioned in the introduction might prove beneficial to address the usually vast differences in the material properties between fiber and coating materials.

\printbibliography

\end{document}